\documentclass[11pt]{amsart}
\usepackage{amssymb,latexsym,amsmath,amsbsy,amsthm,amsxtra,amsgen,graphicx}
\newfont{\msbm}{msbm10 at 11pt}
\newcommand {\R} {\mbox{\msbm R}}
\newcommand {\Z} {\mbox{\msbm Z}}

\newcommand {\N} {\mbox{\msbm N}}

\newfont{\msbmsm}{msbm10 at 8pt}

\newtheorem{Theo}{Theorem}[section]
\newtheorem{Lemma}[Theo]{Lemma}
\newtheorem{Cor}[Theo]{Corollary}

\newtheorem{Rmk}[Theo]{Remark}

\begin{document}
\title[LERW on finite graphs]{Loop-erased random walk on finite graphs and the Rayleigh process}

\author{Jason Schweinsberg}
\address{Department of Mathematics \\ U.C. San Diego \\ 9500 GIlman Drive \\ La Jolla, CA 92093-0112}
\email{jschwein@math.ucsd.edu}
\thanks{Supported in part by NSF Grant DMS-0504882}
\date{July 29, 2007}

\begin{abstract}
Let $(G_n)_{n=1}^{\infty}$ be a sequence of finite graphs, and
let $Y_t$ be the length of a loop-erased random walk on $G_n$ after $t$ steps.  We show that for a large family of sequences of finite graphs, which includes the case in which $G_n$ is the $d$-dimensional torus of size-length $n$ for $d \geq 4$, the process $(Y_t)_{t=0}^{\infty}$, suitably normalized, converges to the Rayleigh process introduced by Evans, Pitman, and Winter.  Our proof relies heavily on ideas of Peres and Revelle, who used loop-erased random walks to show that the uniform spanning tree on large finite graphs converges to the Brownian continuum random tree of Aldous.
\end{abstract}
\maketitle

\section{Introduction}
The loop-erased random walk is a process obtained from a random walk by erasing loops in chronological order.  More precisely, given a sequence of points $\lambda = (u_0, u_1, \dots, u_j)$, which we can think of as the first $j+1$ points visited by some random walk, we define the loop-erasure $LE(\lambda)$ to be the sequence $(v_0, \dots, v_k)$ obtained inductively as follows.  First set $v_0 = u_0$.  Suppose $v_0, \dots, v_m$ have been defined for some $m \geq 0$.  If $v_m = u_j$, then $k = m$ and $v_m$ is the last vertex in the sequence $LE(\lambda)$.  Otherwise, define $v_{m+1} = u_{r+1}$, where $r = \max\{i: u_i = v_m\}$.  We denote the number of points in $LE(\lambda)$, which in this example is $k+1$, by $|LE(\lambda)|$.  We call $|LE(\lambda)|$ the length of the loop-erased path.

The loop-erased random walk was first studied in 1980 by Lawler \cite{law80}, and the model has continued to receive attention in recent years, in part because of connections with uniform spanning trees that were discovered by Pemantle \cite{pem91} and Wilson \cite{wilson}.  For the loop-erased random walk on $\Z^d$ with $d \geq 5$, Lawler \cite{law80} showed that a positive fraction of the vertices never get erased, so if the random walk is run for time $n$, then the length of the loop-erased path is also of order $n$, and the loop-erased random walk as a process converges to Brownian motion.  The loop-erased random walk on $\Z^4$ also converges to Brownian motion, as shown by Lawler in \cite{law86}, but there is a logarithmic correction to the length of the path.  If the random walk is run for time $n$, the length of the loop-erased walk was shown by Lawler \cite{law95} to be of the order $n/(\log n)^{1/3}$.  The loop-erased random walk behaves much differently in dimensions two and three, but there has been recent progress in these lower dimensions.  Kenyon \cite{kenyon} showed that the length of the loop-erased random walk on $\Z^2$ is of order $n^{5/8}$, while Lawler, Schramm, and Werner \cite{lsw} showed that the loop-erased random walk on $\Z^2$ converges to the Schramm-Loewner Evolution (SLE) with parameter $\kappa = 2$.  Kozma \cite{kozma} established the existence of a scaling limit for the loop-erased random walk on $\Z^3$, but the form of the limiting process remains unknown.  In this paper, we consider the behavior of the loop-erased random walk on large finite graphs.  We will focus especially on the $d$-dimensional torus of side length $n$, which we denote by $\Z^d_n = \{0, 1, \dots, n-1\}^d$, for $d \geq 4$.

Given a finite connected graph $G_n = (V_n, E_n)$, write $v \in G_n$ if $v$ is a vertex of $G_n$ and let $|G_n|$ denote the number of vertices of $G_n$.  Write $v \sim w$ if the vertices $v$ and $w$ are connected by an edge.  Throughout the paper, we will assume that $G_n$ is vertex transitive.  Therefore, every vertex has the same degree, which we denote here by $d$.  Let $(X_t)_{t = 0}^{\infty}$ be a discrete-time Markov chain taking its values in $V_n$ such that
\begin{equation}
P(X_{t+1} = w|X_t = v) = \left\{
\begin{array}{ll} 1/2
& \mbox{ if }v = w \\
1/2d &
\mbox{ if }v \sim w \\
0 & \mbox{ otherwise.}
\end{array} \right.
\nonumber
\end{equation}
That is, at each step the Markov chain stays at its current vertex with probability $1/2$ and otherwise moves to a randomly chosen neighboring vertex.  This process is often called the lazy random walk, as opposed to the simple random walk which never stays at its current vertex.  Fix a vertex $o \in G_n$ to be the starting point for the random walk, and denote the transition probabilities of the random walk by $p_{t,n}(x) = P(X_t = x|X_0 = o)$.  Since $G_n$ is vertex transitive, the stationary distribution of the random walk is given by $\pi(x) = 1/|G_n|$ for all $x \in G_n$.  Because the lazy random walk is aperiodic, the distribution of $X_t$ converges to the stationary distribution as $t \rightarrow \infty$.  Denote the uniform mixing time of the random walk, that is, the mixing time measured in terms of the separation distance, by
\begin{equation}\label{mixdef}
\tau_n = \min \bigg\{t: \sup_{x \in G_n} \bigg| \frac{p_{t,n}(x)}{\pi(x)} - 1 \bigg| \leq \frac{1}{2} \bigg\}.
\end{equation}
It is well-known that on $\Z^d_n$, there exist positive constants $C_1$ and $C_2$ such that for all $n$, 
\begin{equation}
C_1 n^2 \leq \tau_n \leq C_2 n^2
\label{taun2}
\end{equation}
(see, for example, the calculations in chapter 5 of \cite{aldfill}).
For the rest of the paper, we will work with a sequence of vertex-transitive, finite graphs $(G_n)_{n=1}^{\infty}$ such that $$\lim_{n \rightarrow \infty} |G_n| = \infty.$$
We will consider two cases.

\bigskip
\noindent {\bf Case 1}:  The graphs $(G_n)_{n=1}^{\infty}$ satisfy the following conditions:
\begin{itemize}
\item There is a constant $C$ such that $$\sup_n \sup_{x \in G_n} \sum_{t=0}^{\lfloor |G_n|^{1/2} \rfloor} (t+1) p_{t,n}(x) \leq C.$$

\item There is a $\delta > 0$ such that
\begin{equation}
\lim_{n \rightarrow \infty} \frac{\tau_n}{|G_n|^{1/2 - \delta}} = 0.
\label{cond2}
\end{equation}
\end{itemize}

\noindent {\bf Case 2}:  For all $n$, we have $G_n = \Z^4_n$.

\bigskip
The conditions for case 1 are precisely the conditions assumed in \cite{perrev} by Peres and Revelle, who used loop-erased random walks to show that the scaling limit of the uniform spanning tree on these graphs as $n \rightarrow \infty$ is the Brownian continuum random tree of Aldous \cite{CRT1}.  As pointed out in \cite{perrev}, this family of graphs includes the $d$-dimensional torus $\Z_n^d$ for $d \geq 5$, the complete graph on $n$ vertices, the hypercubes $\Z_2^n$, and expander graphs.  Schweinsberg \cite{schw} showed that the scaling limit of the uniform spanning tree on $\Z_n^4$ is also the continuum random tree.  The results in this paper will hold for the four-dimensional torus as well as for the family of graphs studied by Peres and Revelle, but at times in the proofs the two cases will be treated separately.

Our goal is to study how the length of the loop-erased random walk on $G_n$ evolves over time.  On $\Z^d$ for $d \geq 5$, it is known that the length of the loop-erased random walk grows linearly in time; see Theorem 7.7.2 of \cite{lawler}, and see also Theorem 7.7.5 of \cite{lawler} for a similar result when $d = 4$.  The reason is that the random walk does not make long loops, so loop erasure is a local procedure.  The linear growth then comes from the fact that most points have approximately the same probability of never being erased.  For the tori $\Z^d_n$ with $d \geq 4$ and for the other graphs satisfying the conditions of case 1, the random walk also makes long loops.  The long loops occur on a time scale much longer than the mixing time.  As a result, when the random walk makes a long loop, the point on the path that the random walk hits is approximately uniformly distributed over all the points in the path.  Therefore, when the long loop is erased, the length of the path is multipied by a fraction which is approximately uniformly distributed between $0$ and $1$.  As on $\Z^d$, the length of the path grows approximately linearly between the times when these long loops form.

As a result of this intuition, Jim Pitman conjectured that the length of the loop-erased random walk on $G_n$ converges to a process called the Rayleigh process as $n \rightarrow \infty$.  The Rayleigh process $(R(t), t \geq 0)$ was introduced by Evans, Pitman, and Winter \cite{evpitwin}.  The process grows linearly at unit speed between jumps.  At time $t$, jumps occur at rate $R(t-)$, and at the times of jumps, the value of the process gets multiplied by a random variable which is uniformly distributed on $[0,1]$.  More formally, the Rayleigh process can be constructed from a Poisson point process $\Pi$ on $[0, \infty) \times [0, \infty)$ whose intensity measure is Lebesgue measure.  Given $y \geq 0$, we can obtain a Rayleigh process started from $R(0) = y$ by defining
\begin{equation}
R(t) = (y + t) \wedge \inf\{x + (t-s): (s, x) \in \Pi, 0 \leq s \leq t\}.
\label{Rteq}
\end{equation}
This means that when $(t,x)$ is a point of this Poisson process and $R(t-) > x$, there is a jump at time $t$ and $R(t) = x$ (see Figure 1 below).  It was shown in Proposition 8.1 of \cite{evpitwin} that the stationary distribution for this process is the Rayleigh distribution, where we say a random variable $W$ has a Rayleigh distribution if $P(W > x) = e^{-x^2/2}$ for all $x \geq 0$, and that for any starting point $y$, the distribution of $R(t)$ converges to the Rayleigh distribution as $t \rightarrow \infty$.  For more about the Rayleigh process, see section 8 of \cite{evpitwin}.

It is already known (see \cite{perrev} for case 1 and \cite{schw} for case 2) that if $x$ and $y$ are vertices of $G_n$ chosen uniformly at random, then the distribution of the length of the loop-erased random walk started at $x$ and run until it hits $y$, suitably normalized, converges to the Rayleigh distribution as $n \rightarrow \infty$.  The theorem below, which was conjectured by Pitman, is a dynamical result, which shows that the length of the loop-erased random walk converges to the Rayleigh process, in the sense of Skorohod convergence for processes whose sample paths are right continuous and have left limits.  When $G_n$ is the complete graph with $n$ vertices, this result can be obtained from Corollary 8.2 of \cite{evpitwin}.

\begin{figure}
\includegraphics{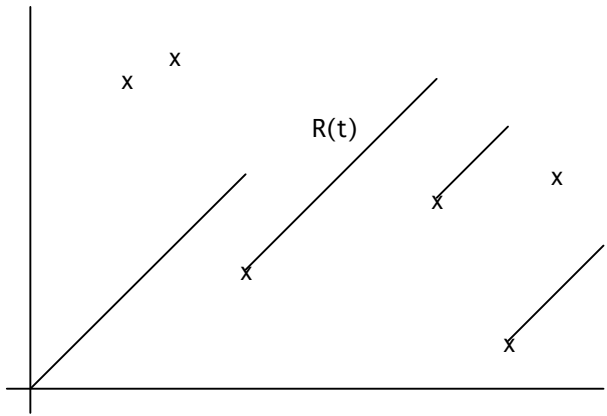}
\caption{Figure 1: The Rayleigh process}
\vspace{-8.8in}
Figure 1: The Rayleigh process
\end{figure}

\begin{Theo}
Let $R = (R(t), t \geq 0)$ denote the Rayleigh process with $R(0) = 0$.
Let $(G_n)_{n=1}^{\infty}$ be a sequence of graphs satisfying the conditions of either case 1 or case 2.  Let $(X_t)_{t=0}^{\infty}$ be a lazy random walk on $G_n$, as defined above.  For all $t$, let $Y_t$ be the length of the loop-erased path $LE((X_s)_{s=0}^t)$.  Then there exist sequences of constants $(a_n)_{n=1}^{\infty}$ and $(b_n)_{n=1}^{\infty}$ satisfying $0 < \inf a_n \leq \sup a_n < \infty$ and $0 < \inf b_n \leq \sup b_n < \infty$ such that if we
define the continuous-time process $Z_n = (Z_n(t), t \geq 0)$ by
\begin{equation}
Z_n(t) = \left\{
\begin{array}{ll} 
b_n |G_n|^{-1/2} Y_{\lfloor a_n |G_n|^{1/2} t \rfloor} & \mbox{ in case 1 } \\
b_n n^{-2} (\log n)^{-1/6} Y_{\lfloor a_n n^2 (\log n)^{1/2} t \rfloor} & \mbox{ in case 2,}
\end{array} \right.
\nonumber
\end{equation}
then $Z_n \rightarrow R$ in the Skorohod topology as $n \rightarrow \infty$.
\label{Rayleighthm}
\end{Theo}

\begin{Rmk}
{\em We work with the lazy random walk rather than the simple random walk because our proof uses results about the mixing time which require the random walk to be aperiodic.  However, once convergence to the Rayleigh process is established for the lazy random walk, it follows easily for the simple random walk with $a_n$ replaced by $a_n/2$.}
\end{Rmk}

\begin{Rmk}
{\em By the results in section 8 of \cite{perrev}, when $G_n = \Z^d_n$ for some $d \geq 5$, there are positive numbers $a$ and $b$ such that $\lim_{n \rightarrow \infty} a_n = a$ and $\lim_{n \rightarrow \infty} b_n = b$.  It follows from the definition of the Skorohod metric that we still have $Z_n \rightarrow R$ in the Skorohod topology if we take $a_n = a$ and $b_n = b$ for all $n$.}
\end{Rmk}

\begin{Rmk}
{\em To understand the scaling in case 2, note that the appropriate time scale for convergence to the Rayleigh process is the time scale on which long loops occur, as these long loops correspond to the jumps of the Rayleigh process.  Given a random walk segment of length $L$, if one removes a segment of length $\tau_n$ from the middle, then the probability that the two remaining segments intersect (which would correspond to a long loop) is the same order of magnitude as the probability that two independent random walk segments of length $L$ on $\Z^4_n$ started from the uniform distribution intersect.  The expected number of intersection times between two such walks is $L^2/n^4$.  Furthermore, if there is one intersection, then there will likely be $O(\log n)$ intersections because two independent random walks of length $n$ in $\Z^4$ started at the origin will intersect $O(\log n)$ times (see Proposition 3.2.3 of \cite{lawler}).  Consequently, the intersection probability of two random walks of length $L$ on $\Z^4_n$ started from the uniform distribution is of order $L^2/(n^4 \log n)$ (see Propositions 2.8 and 2.9 of \cite{schw}), which is of order one when $L$ is of order $n^2 (\log n)^{1/2}$.  This explains the time scaling in case 2.  The spatial scaling results from the fact that the number of points remaining after loop-erasure is of order $n^2 (\log n)^{1/2}/(\log n)^{1/3} = n^2 (\log n)^{1/6}$.
}
\end{Rmk}

There are two steps to proving Theorem \ref{Rayleighthm}.  First, we must show the convergence of finite-dimensional distributions.  That is, we need to show that if $0 \leq t_1 < t_2 < \dots < t_k$, then $(Z_n(t_1), \dots, Z_n(t_k))$ converges weakly to $(R(t_1), \dots, R(t_k))$.  This is done in section \ref{fddsec}.  Then, we must show that the sequence of processes $(Z_n)_{n=1}^{\infty}$ is relatively compact, which we do in section \ref{tightsec}.  These results imply that $Z_n \rightarrow R$ in the Skorohod topology (see Theorem 7.8 in chapter 3 of \cite{ek}).  

\section{Convergence of finite-dimensional distributions}
\label{fddsec}

Fix $k \in \N$, and fix times $0 \leq t_1 < t_2 < \dots < t_k$.  Our goal in this section is to prove that
\begin{equation}
(Z_n(t_1), \dots, Z_n(t_k)) \rightarrow_d (R(t_1), \dots, R(t_k)).
\label{fddeq}
\end{equation}
The proof proceeds in three steps.  First, we review some results concerning the loop-erased random walk on $G_n$.  Next we set up a coupling between the loop-erased random walk and the Rayleigh process.  The result (\ref{fddeq}) will then follow from some bounds for the Rayleigh process.

\subsection{Loop-erased random walk on $G_n$}

To study the loop-erased random walk on $G_n$, we follow the strategy introduced by Peres and Revelle in \cite{perrev} of splitting the random walk into shorter segments.  The ideas are the same in cases 1 and 2, but it is necessary treat the two cases separately.  Some of the minor differences in the treatments of the two cases could be avoided, but we prefer to set up the notation so that we can directly invoke results in \cite{perrev} and \cite{schw}.

First, consider case 1.  Following \cite{perrev}, let $r = \lfloor \tau_n^{1/4} |G_n|^{3/8} \rfloor$ and $s = \lfloor \tau_n^{3/4} |G_n|^{1/8} \rfloor$.  Note that (\ref{cond2}) implies that $s$ is much smaller than $r$ when $n$ is large.  We will work with segments of the random walk whose length is approximately $r$.  For all positive integers $i$, let $A_i = \{(i-1)r + 2s + 1, (i-1)r + 2s + 2, \dots, ir - s\}$, which will be the set of times associated with the $i$th segment.  Say that a time $u$ is locally retained if $LE((X_t)_{t=\max\{0, u-s\}}^u) \cap (X_t)_{t=u+1}^{u+s} = \emptyset$.  Let $U$ denote the set of times $u \in A_i$ that are locally retained.  Still following \cite{perrev}, define the local loop erasure of the segment $A_i$ by $LE_s(A_i) = (X_t)_{t \in A_i \cap U}$.  Denote by $|LE_s(A_i)|$ the cardinality of $A_i \cap U$, which is the length of the path $LE_s(A_i)$.  For any $V \subset G_n$, define the capacity of the set $V$ by letting $(W_t)_{t=0}^{\infty}$ be a random walk on $G_n$ started from the stationary distribution $\pi$ and then defining
\begin{equation}
\mbox{Cap}(V) = P((W_t)_{t=0}^r \cap V \neq \emptyset).
\label{capdef}
\end{equation}
As in \cite{perrev}, define the constants $\gamma_n = r^{-1} E[|LE_s(A_i)|]$ and $\alpha_n = r^{-2} |G_n| E[\mbox{Cap}(LE_s(A_i))]$.  As can be seen from Lemma 5.3 of \cite{perrev}, the sequences of constants $(\alpha_n)_{n=1}^{\infty}$ and $(\gamma_n)_{n=1}^{\infty}$ are bounded away from zero and infinity.

For case 2, let $r = \lfloor n^2 (\log n)^{9/22} \rfloor$ and let $A_i = \{(i-1)r, (i-1)r + 1, \dots, ir-1\}$ for all positive integers $i$, as in \cite{schw}.  By Corollary 3.2 of \cite{schw}, we can define a sequence of constants $(\gamma_n)_{n=1}^{\infty}$, bounded away from zero and infinity, such that for some positive constants $C$ and $C'$, we have
\begin{equation}
P \bigg( \bigg| \big| LE((X_t)_{t \in A_i}) \big| - \gamma_n n^2 (\log n)^{5/66} \bigg| > \frac{C' n^2}{(\log n)^{15/44}} \bigg) \leq \frac{C}{(\log n)^{2/11}}
\label{lenbound}
\end{equation}
for all $n$.  The exponent of $5/66 = 9/22 - 1/3$ comes from the fact that in four dimensions, the length of the loop erasure of a random walk of length $r$ is of the order $r/(\log r)^{1/3}$.   Also, fix a small positive number $\eta > 0$, and let $w = \lfloor n^2 (\log n)^{\eta} \rfloor.$  For $V \subset \Z^4_n$, define $\mbox{Cap}(V)$ as in (\ref{capdef}) but with $r - 2w$ in place of $r$.  By Proposition 3.8 of \cite{schw}, there is a sequence of constants $(\alpha_n)_{n=1}^{\infty}$, bounded away from zero and infinity, such that for some positive constant $C$, we have
$$P \bigg( \bigg| \mbox{Cap} \big( LE((X_t)_{t \in A_i}) \big) - \frac{\alpha_n}{(\log n)^{2/11}} \bigg| > \frac{1}{(\log n)^{5/22}} \bigg) \leq \frac{C}{(\log n)^{3/22 - \eta}}.$$
These bounds show that the length and capacity of the loop-erased segment $LE((X_t)_{t \in A_i})$ are highly concentrated around their means.

In the absence of long loops, the $i$th segment of the random walk of length $r$, after loop erasure, looks approximately like the path $LE_s(A_i)$ in case 1 and like the path $LE((X_t)_{t \in A_i})$ in case 2.  However, long loops can cause entire segments of length $r$ to get erased.
We will use indicator random variables to keep track of the long loops.  In case 1, for $i < j$, let $I_{i,j}$ be the indicator of the event that $LE_s(A_i) \cap (X_t)_{t \in A_j} \neq \emptyset$.  In case 2, let $I_{i, i+1} = 0$ for all $i$ and, for $i < j - 1$, let $I_{i,j}$ be the indicator of the event that $LE((X_t)_{t \in A_i}) \cap (X_t)_{t \in A_j} \neq \emptyset$.  Thus, in both cases, when $I_{i,j} = 1$, the $j$th segment of the random walk of length $r$ intersects the loop-erasure of the $i$th segment.  For both cases, let $S_0 = \{0\}$ and, for $j \geq 1$, let $$S_j = \{k \in S_{j-1}: I_{i,j} = 0 \mbox{ for all }i \in \{1, \dots, k\} \cap S_{j-1}\} \cup \{j\}.$$  Here $S_j$ consists of the indices of the segments that are not erased in the loop erasure through time $jr$, with the convention that if segment $j$ intersects the loop-erasure of segment $i$, causing segments $i+1, \dots, j-1$ and parts of segments $i$ and $j$ to be erased, we keep $j$ in $S_j$ but not $i$.

The number of segments not erased through time $jr$ is $|S_j|$, where $|S_j|$ denotes the cardinality of $S_j$.  Since $|LE_s(A_i)|$ has mean $\gamma_n r$ in case 1 and $|LE((X_t)_{t \in A_i})|$ has mean approximately $\gamma_n n^2 (\log n)^{5/66}$ in case 2, the length $Y_{jr}$ of the loop-erased random walk at time $jr$ can be approximated by $\gamma_n r |S_j|$ in case 1 and by $\gamma_n n^2 (\log n)^{5/66} |S_j|$ in case 2.  More precisely, we have the following result.

\begin{Lemma}
Let $(T_n)_{n=1}^{\infty}$ be a sequence of fixed times such that for some constants $C_1$ and $C_2$, we have $T_n \leq C_1 |G_n|^{1/2}$ for all $n$ in case 1 or $T_n \leq C_2 n^2 (\log n)^{1/2}$ for all $n$ in case 2.  Let $b_n = \alpha_n^{1/2} \gamma_n^{-1}$.  Then, there are positive constants $C$ and $C'$ such that for all $n$, 
\begin{align}
&P \bigg( \bigg| \frac{Y_{T_n} - \gamma_n r |S_{\lceil T_n/r \rceil}|}{b_n^{-1} |G_n|^{1/2}} \bigg| > |G_n|^{-\delta/24} \bigg) \leq C |G_n|^{-3 \delta/16} &\mbox{in case 1,} \nonumber \\
&P \bigg( \bigg| \frac{Y_{T_n} - \gamma_n n^2 (\log n)^{5/66} |S_{\lceil T_n/r \rceil}|}{b_n^{-1} n^2 (\log n)^{1/6}} \bigg| > C'(\log n)^{-1/22} \bigg) \leq \frac{C (\log \log n)^2}{(\log n)^{1/22}} &\mbox{in case 2.} \nonumber
\end{align}
\label{prlem}
\end{Lemma}

\begin{proof}
Equation (41) in section 6 of \cite{perrev} gives the estimate for case 1.  Although (41) in \cite{perrev} is stated for certain random times rather than fixed times, the arguments leading to this result (which show, for example, that $|LE_s(A_i)|$ is highly concentrated around its mean, and that the contributions of the gaps of length $3s$ between the segments $A_i$ can be neglected) hold for fixed $T$ as well.

The result for case 2 follows from (\ref{lenbound}), and from Proposition 3.10, Proposition 4.1, and Lemma 4.11 of \cite{schw}.  Proposition 3.10 of \cite{schw} shows that on the complement of an event whose probability is shown in Proposition 4.1 of \cite{schw} to be at most $C(\log n)^{-1/11}$, the only discrepancies between $|LE((X_t)_{t=0}^{T_N})| = Y_{T_N}$ and $\gamma_n n^2 (\log n)^{5/66} |S_{\lceil T_n/r \rceil}|$ come from the following three sources:
\begin{itemize}
\item There are differences, for $i \in S_{\lceil T_n/r \rceil}$, between $|LE((X_t)_{t \in A_i})|$ and $\gamma_n n^2 (\log n)^{5/66}$.  By (\ref{lenbound}), on the complement of an event of probability at most $C \lceil T_n/r \rceil (\log n)^{-2/11}$, which is of order $(\log n)^{-1/11}$, the sum of the absolute values of these differences can be bounded above by $C' \lceil T_n/r \rceil n^2 (\log n)^{-15/44}$, which is of order $n^2 (\log n)^{-15/44 + 1/11} = n^2 (\log n)^{-1/4}$.

\item There are differences, for $i \in S_{\lceil T_n/r \rceil}$, between $|LE((X_t)_{t \in A_i})|$ and the number of points from the $i$th segment of the random walk that end up in the path $LE((X_t)_{t=0}^{T_n})$.  By Proposition 3.10 of \cite{schw}, these differences add up to at most $2 |S_{\lceil T_n/r \rceil}| w$, which is of order $n^2 (\log n)^{1/11 + \eta}$.

\item There are contributions from segments that get only partially erased because of intersections.  By Lemma 4.11 of \cite{schw}, on the complement of an event whose probability is of order $(\log \log n)^2 (\log n)^{-1/22}$, there are at most $(\log n)^{1/22}$ such contributions, each of order $n^2 (\log n)^{5/66}$, so the total contribution is of order $n^2 (\log n)^{1/22 + 5/66} = n^2 (\log n)^{1/6 - 1/22}$.
\end{itemize}
These observations give the result for case 2.
\end{proof}

Let $a_n = \alpha_n^{-1/2}$, and let $b_n = \alpha_n^{1/2} \gamma_n^{-1}$ as in Lemma \ref{prlem}.  
Let
\begin{equation}
d_n = \left\{
\begin{array}{ll} 
r \alpha_n^{1/2} |G_n|^{-1/2} & \mbox{ in case 1 } \\
\alpha_n^{1/2} (\log n)^{-1/11} & \mbox{ in case 2.}
\end{array} \right.
\nonumber
\end{equation}
Here $d_n$ is an estimate of the length of an individual loop-erased segment after the rescaling in Theorem \ref{Rayleighthm}.
Recall that our goal is to find the limit in distribution of $(Z_n(t_1), \dots, Z_n(t_k))$.  The next corollary, which follows immediately from Lemma \ref{prlem}, relates $Z_n(t_i)$ to the sets $S_j$ by showing that $Z_n(t_i)$ can be approximated by the product of $|S_{t_i^*}|$ (the number of segments retained) and $d_n$ (the length of a segment).

\begin{Cor}\label{lerwlem}
For $i = 1, \dots, k$, let $t_i^* = \lceil T_n/r \rceil$, where $T_n = \lfloor a_n |G_n|^{1/2} t_i \rfloor$ in case 1 and $T_n = \lfloor a_n n^2 (\log n)^{1/2} t_i \rfloor$ in case 2.
There exist positive constants $C$ and $C'$, depending on $k$, such that for all $i = 1, \dots, k$ and all $n$, we have
\begin{align}
&P \big( \big| Z_n(t_i) - d_n |S_{t_i^*}| \big| > |G_n|^{-\delta/24} \big) \leq C |G_n|^{-3 \delta/16} &\mbox{in case 1,} \nonumber \\
&P \big( \big| Z_n(t_i) - d_n |S_{t_i^*}| \big| > C' (\log n)^{-1/22} \big) \leq \frac{C (\log \log n)^2}{(\log n)^{1/22}} &\mbox{in case 2.} \nonumber
\end{align}
In particular, in both cases,
\begin{equation}\label{convprob1}
(Z_n(t_1) - d_n |S_{t_1^*}|, \dots, Z_n(t_k) - d_n |S_{t_k^*}|) \rightarrow_p 0,
\end{equation}
where $\rightarrow_p$ denotes convergence in probability, using the Euclidean metric on $\R^k$.
\end{Cor}

\subsection{Coupling with the Rayleigh process}

We now show how the loop-erased random walk on $G_n$ can be coupled with the Rayleigh process.  It will suffice to couple the sets $S_j$.  We begin by reviewing the coupling between the loop-erased random walk on $G_n$ and loop-erased random walk on the complete graph used in \cite{perrev} and \cite{schw}.  Let $m = \lceil r^{-2} \alpha_n^{-1} |G_n| \rceil = \lceil d_n^{-2} \rceil$ in case 1, and let $m = \lfloor \alpha_n^{-1} (\log n)^{2/11} \rfloor = \lfloor d_n^{-2} \rfloor$ in case 2.  Let $(\xi_i)_{i=1}^{\infty}$ be an i.i.d. sequence of random variables which have the uniform distribution on $\{1, \dots, m\}$.  The process $(\xi_i)_{i=1}^{\infty}$ can be viewed as a lazy version of a random walk on the complete graph with $m$ vertices in which on each step the random walk stays in its current position with probability $1/m$.  For $i < j$, let ${\tilde I}_{i,j} = {\bf 1}_{\{\xi_i = \xi_j\}}$.  Let ${\tilde S}_0 = \{0\}$ and, for $j \geq 1$, let $${\tilde S}_j = \{k \in {\tilde S}_{j-1}: {\tilde I}_{i,j} = 0 \mbox{ for all }i \in \{1, \dots, k\} \cap {\tilde S}_{j-1}\} \cup \{ j \}.$$  Here ${\tilde S}_j$ consists of the vertices of the walk that are not erased after $j$ steps.

The lemma below shows that the sets $S_j$ can be coupled with the sets ${\tilde S}_j$ with high probability.  This coupling is possible because the $j$th segment $(X_t)_{t \in A_j}$ is approximately equally likely to intersect any of the previous segments.  This is true because the mixing time $\tau_n$ is much shorter than the length of the segments (a consequence of (\ref{cond2}) in case 1 and (\ref{taun2}) in case 2), so conditional on the previous segments, the $j$th segment is at approximately a uniform random point a short distance into the segments.  The proof of Lemma 6.3 of \cite{perrev} gives the bound for case 1, while the bound for case 2 comes from Proposition 4.7 of \cite{schw}.  The results in \cite{perrev} and \cite{schw} are stated for the case in which the random walk is run for a random time rather than the fixed time $t_k^*$, but the same proofs work for fixed times.

\begin{Lemma}
Define $t_k^*$ as in Corollary \ref{lerwlem}.  There exists a coupling of $(S_j)_{j=1}^{\infty}$ and $({\tilde S}_j)_{j=1}^{\infty}$ such that for some constant $C$, we have
\begin{equation}
P(S_j = {\tilde S}_j \mbox{ for all }j = 1, \dots, t_k^*) \geq \left\{
\begin{array}{ll} 
1 - C|G_n|^{-\delta/16} & \mbox{ in case 1 } \\
1 - C (\log n)^{-1/22 + \eta} & \mbox{ in case 2.}
\end{array} \right.
\nonumber
\end{equation}
\label{prcouple}
\end{Lemma}

Our next step is to couple the ${\tilde S}_j$ with a Rayleigh process.  Let $\Pi$ be a Poisson point process with Lebesgue intensity on $[0, \infty) \times [0, \infty)$.  For all $t \geq 0$, define $R_t$ by (\ref{Rteq}) with $y = 0$.  For $0 \leq i < j$, let $I'_{i,j}$ be the indicator of the event that there is at least one point of $\Pi$ in
\begin{equation}
[d_n (j-1), d_n j) \times [d_n (i-1), d_n i)
\label{rect}
\end{equation}
Let $S'_0 = \{0\}$.  For $j \geq 1$, conditional on $S_{j-1}' = \{\ell_1, \dots, \ell_{|S_{j-1}'|}\}$, where $\ell_1 < \dots < \ell_{|S_{j-1}'|}$, define
$$S'_j = \{\ell_k \in S'_{j-1}: I'_{i,j} = 0 \mbox{ for all }i \leq k\} \cup \{ j \}.$$

Because the rectangles in (\ref{rect}) have area $d_n^2$ and are disjoint, the indicator random variables $I'_{i,j}$ are independent and equal $1$ with probability $1 - e^{-d_n^2}$.  Also, conditional on ${\tilde S}_{j-1} = \{\ell_1, \dots, \ell_{|{\tilde S}_{j-1}|}\}$, the probability that ${\tilde I}_{i,j} = 0$ for all $i \in {\tilde S}_{j-1}$ is $1 - |{\tilde S}_{j-1}|/m$ and, for all $h \leq |{\tilde S}_{j-1}|$, the probability that ${\tilde I}_{\ell_h, j} = 1$ and $\tilde{I}_{i,j} = 0$ for all $i \in {\tilde S}_{j-1}$ with $i \neq h$ is $1/m$.  Consequently, we will be able to couple the ${\tilde S}_j$ and $S'_j$ by using the following elementary lemma.  Because this result is a special case of Lemma 4.5 in \cite{schw}, we omit the proof.
 
\begin{Lemma}
Suppose $0 < p < 1/j$ and $0 < q < 1$.  Suppose $V_1, \dots, V_j$ are random variables such that $P(V_i = 0 \mbox{ for all }i) = 1 - jp$ and $P(V_i = 1 \mbox{ and }V_{\ell} = 0 \mbox{ for }\ell \neq i) = p$ for $i = 1, \dots, j$.  Suppose $W_1, \dots, W_j$ are independent random variables such that $P(W_i = 1) = q$ and $P(W_i = 0) = 1-q$.  Then there is a coupling of $V_1, \dots, V_j$ and $W_1, \dots, W_j$ such that
$$P(V_i = W_i \mbox{ for all }i) \geq 1 - j|p-q| - j(j-1)q^2.$$
\label{easycouple}
\end{Lemma}

\begin{Lemma}
Define $t_k^*$ as in Corollary \ref{lerwlem}. 
There exists a coupling of $(S'_j)_{j=1}^{\infty}$ and $({\tilde S}_j)_{j=1}^{\infty}$ such that, for some constant $C$, we have 
\begin{equation}
P(S'_j = {\tilde S}_j \mbox{ for all }j = 1, \dots, t_k^*) \geq \left\{
\begin{array}{ll} 
1 - C|G_n|^{-\delta/4} & \mbox{ in case 1 } \\
1 - C (\log n)^{-1/11} & \mbox{ in case 2.}
\end{array} \right.
\nonumber
\end{equation}
\label{raycouple}
\end{Lemma}

\begin{proof}
Let $p = 1/m$ and $q = 1 - e^{-d_n^2}$.  Note that $|S'_{j-1}| \leq j$ and $|{\tilde S}_{j-1}| \leq j$ for all $j$.  By Lemma \ref{easycouple}, conditional on $S'_{j-1} = {\tilde S}_{j-1}$, a coupling can be achieved such that the probability that $S'_j = {\tilde S}_j$ is at least $1 - j|p-q| - j(j-1)q^2$.  Therefore, there is a coupling such that
\begin{align}
P(S'_j = {\tilde S}_j \mbox{ for all }j = 1, \dots, t_k^*) &\geq 1 - \sum_{j=1}^{t_k^*} \big( j|p-q| - j(j-1)q^2 \big) \nonumber \\
&\geq 1 - (t_k^*)^2 |p-q| - (t_k^*)^3 q^2.
\label{couplebound}
\end{align}
For all $x > 0$, we have $0 \leq e^{-x} - 1 + x \leq x^2/2$ and $0 \leq 1/x - 1/\lceil x \rceil \leq 1/x - 1/(x+1) \leq 1/x^2$.  Also, if $x \geq 2$, then $0 \leq 1/\lfloor x \rfloor - 1/x \leq 1/(x-1) - 1/x \leq 2/x^2$.  Therefore, we have in case 1
$$|p - q| = \bigg| \frac{1}{\lceil d_n^{-2} \rceil} - (1 - e^{-d_n^2}) \bigg| \leq |e^{-d_n^2} - 1 + d_n^2| + \bigg| d_n^2 - \frac{1}{\lceil d_n^{-2} \rceil} \bigg| \leq \frac{d_n^4}{2} + d_n^4 = \frac{3d_n^4}{2}$$
and in case 2 for $n$ large enough that $d_n^{-2} \geq 2$,
$$|p - q| = \bigg| \frac{1}{\lfloor d_n^{-2} \rfloor} - (1 - e^{-d_n^2}) \bigg| \leq |e^{-d_n^2} - 1 + d_n^2| + \bigg| d_n^2 - \frac{1}{\lfloor d_n^{-2} \rfloor} \bigg| \leq \frac{d_n^4}{2} + 2 d_n^4 = \frac{5d_n^4}{2}.$$
Also, $q^2 = (1 - e^{-d_n^2})^2 \leq d_n^4.$  By combining these bounds with (\ref{couplebound}), we can bound the probability that $S_j' \neq {\tilde S}_j$ for some $j \leq t_k^*$ in both cases by
$$\frac{5}{2} (t_k^*)^2 d_n^4 + (t_k^*)^3 d_n^4.$$  Because there is a constant $C$ such that $t_k^* \leq C r^{-1} |G_n|^{1/2}$ in case 1 and $t_k^* \leq C (\log n)^{1/11}$ in case 2, the result follows.
\end{proof}

\begin{Rmk}\label{refrem}
{\em As another way of understanding the scaling, let $L_j$ be the cardinality of ${\tilde S}_j$, which is the length of the loop-erased random walk on the complete graph $K_m$ after $j$ steps.  Conditional on $L_j$, with probability $1 - L_j/m$ the next step of the walk will not form a loop and we will have $L_{j+1} = L_j + 1$, and for $k = 1, \dots, L_j$, with probability $1/m$ the next step of the walk will duplicate the $k$th vertex on the current path, and after the loop erasure we will have $L_j = k$.  Thus,
$$E[L_{j+1}|L_j] = \bigg(1 - \frac{L_j}{m} \bigg)(L_j + 1) + \frac{1}{m} \sum_{k=1}^{L_j} k = 1 + L_j - \frac{L_j}{2m} - \frac{L_j^2}{2m},$$ and so $E[L_{j+1}] = 1 + E[L_j] - E[L_j]/2m - E[L_j^2]/2m$.
Therefore, letting $x = \lim_{j \rightarrow \infty} E[L_j]$, we have $x \approx 1 + x - x/2m - x^2/2m$, and since $x^2$ is much larger than $x$, it follows that $x \approx \sqrt{2m}$, where the approximations only give the correct order of magnitude because $x^2$ is being used to approximate $E[L_j^2]$.  Thus, in the long run, the loop-erased paths should contain on the order of $\sqrt{m}$ segments, which in case 1 is the order of $d_n^{-1}$ and in case 2 is of order $(\log n)^{1/11}$.  Loop-erased segments have length of order $r$ in case 1 and of order $n^2 (\log n)^{5/66}$ in case 2, so the length of the loop-erased path should be of order $r/d_n$, which is of order $|G_n|^{1/2}$, in case 1 and of order $n^2 (\log n)^{5/66} (\log n)^{1/11} = (\log n)^{1/6}$ in case 2, consistent with the scaling in Theorem \ref{Rayleighthm}.}
\end{Rmk}

\subsection{Bounds for the Rayleigh process}

Corollary \ref{lerwlem} bounds the process $Z_n$ using the sets $S_j$, and Lemmas \ref{prcouple} and \ref{raycouple} couple the sets $S_j$ and the sets $S_j'$.  In this subsection, we complete the proof of (\ref{fddeq}) by using the sets $S_j'$ to obtain bounds for the Rayleigh process $R$.  Recall that both the Rayleigh process and the sets $S_j'$ are constructed from the same Poisson process $\Pi$.  We begin by obtaining a bound which is valid at times that are integer multiples of $d_n$.  This is a deterministic bound which follows from the construction of $R$.

\begin{Lemma}
For all nonnegative integers $j$, we have $\big| R(d_n j) - d_n |S_j'| \big| \leq d_n$.
\label{detraylem}
\end{Lemma}

\begin{proof}
We proceed by induction on $j$.  The result is trivial for $j = 0$.  Assume, for some integer $j \geq 1$, the result holds for $j-1$.  We consider two cases.

First, suppose $I_{i,j}' = 0$ for all $i \leq |S_{j-1}'|$.  Then $|S_j'| = |S_{j-1}'| + 1$.  By the definition of the $I_{i,j}'$, there is no point of $\Pi$ in $[d_n (j-1), d_n j) \times [0, d_n |S_{j-1}'|)$.  It follows that $$\min \big\{ d_n |S_{j-1}'|, R(d_n(j-1)) + d_n \big\} \leq R(d_n j) \leq R(d_n (j-1)) + d_n.$$  Now the induction hypothesis gives $d_n(|S_j'| - 1) \leq R(d_n j) \leq d_n(|S_j'| - 1) + 2d_n,$ which implies the claim.

Alternatively, suppose there is an $\ell \leq |S_{j-1}'|$ such that  $I_{\ell, j}' = 1$ and $I_{i,j}' = 0$ for $i < \ell$.  Then $|S'_j| = \ell$.  There is no point of $\Pi$ in $[ d_n (j-1), d_n j) \times [0, d_n (\ell - 1))$ but there is a point of $\Pi$ in $[ d_n (j-1), d_n j) \times [d_n (\ell - 1), d_n \ell)$.  Therefore, one can see from the construction that $$d_n (\ell - 1) \leq R(d_n j) \leq d_n \ell + d_n,$$
and again the claim follows.
\end{proof}

\begin{Lemma}
There is a constant $C$ such that for all $n$ and for $i = 1, \dots, k$, we have
\begin{equation}
P \big( \big| R(t_i) - d_n |S'_{t_i^*}| \big| > 3d_n \big) \leq \left\{
\begin{array}{ll} 
C|G_n|^{-\delta/4} & \mbox{ in case 1 } \\
C (\log n)^{-1/11} & \mbox{ in case 2.}
\end{array} \right.
\nonumber
\end{equation}
In particular
\begin{equation}\label{convprob2}
(d_n |S_{t_1^*}'|, \dots, d_n |S_{t_k^*}'|) \rightarrow_p (R(t_1), \dots, R(t_k)).
\end{equation}
\label{mainraylem}
\end{Lemma}

\begin{proof}
From the definitions of $t_i^*$ and $d_n$, we get $d_n (t_i^* - 2) \leq t_i \leq d_n (t_i^* + 1)$ in both cases for large enough $n$.  For all $0 \leq s < t$, we have $R(t) \leq R(s) + (t-s)$.  Using this fact for the first inequality and Lemma \ref{detraylem} for the second, we get
$$R(t_i) \leq R(d_n (t_i^* - 2)) + 3 d_n \leq d_n |S'_{t_i^* - 2}| + 4 d_n.$$
Likewise, this time using Lemma \ref{detraylem} for the first inequality,
$$d_n |S'_{t_i^* + 1}| - 4 d_n \leq R(d_n (t_i^* + 1)) - 3 d_n \leq R(t_i).$$
Therefore, on the event that
\begin{equation}
|S'_{t_i^* - 2}| + 2 = |S'_{t_i^*}| = |S'_{t_i^* + 1}| - 1,
\label{nojumpevent}
\end{equation}
we have $$\big| R(t_i) - d_n |S'_{t_i^*}| \big| \leq 3 d_n.$$  Thus, it remains only to bound the probability that (\ref{nojumpevent}) fails to occur.  However, the event (\ref{nojumpevent}) occurs as long as $I'_{\ell, j} = 0$ whenever $j \in \{t_i^*-1, t_i^*, t_i^* + 1\}$ and $\ell \leq t_i^* + 1$.  Recall that the random variables $I'_{\ell, j}$ are nonzero with probability $1 - e^{-d_n^2} \leq d_n^2$, so the probability that (\ref{nojumpevent}) fails to hold is at most $3 (t_i^* + 1) d_n^2$.  The bounds in the lemma now from the definition of $d_n$ and the fact that there is a constant $C$ such that $t_k^* \leq C r^{-1} |G_n|^{1/2}$ in case 1 and $t_k^* \leq C (\log n)^{1/11}$ in case 2, and the convergence in probability follows easily from these bounds.
\end{proof}
 
\begin{proof}[Proof of (\ref{fddeq})]
Recall that if $(S, d)$ is a metric space and $\mu$ and $\nu$ are probability measures on the Borel $\sigma$-field of $(S, d)$, then the Prohorov distance between $\mu$ and $\nu$ is defined by $$\rho(\mu, \nu) = \inf\{\epsilon > 0: \mu(F) \leq \nu(F^{\epsilon}) + \epsilon \mbox{ for all closed subsets }F \subset S\},$$ where $F^{\epsilon} = \{x \in S: d(x, y) < \epsilon \mbox{ for some }y \in S\}$.  It is well-known that if $X, X_1, X_2, \dots$ are $S$-valued random variables with distributions $\mu, \mu_1, \mu_2, \dots$ respectively, then $X_n \rightarrow_d X$ as $n \rightarrow \infty$ if and only if $\rho(\mu_n, \mu) \rightarrow 0$ as $n \rightarrow \infty$ (see, for example, Theorem 3.1 in chapter 3 of \cite{ek}).

Now consider the metric space $\R^k$ with the Euclidean metric.  Let $\mu_n$ and $\nu_n$ denote the distributions of $(d_n |S_{t_1^*}|, \dots, d_n |S_{t_k^*}|)$ and $(d_n |S_{t_1^*}'|, \dots, d_n |S_{t_k^*}'|)$ respectively.  Lemmas \ref{prcouple} and \ref{raycouple} give that $\rho(\mu_n, \nu_n) \rightarrow 0$ as $n \rightarrow \infty$.  Convergence in probability implies convergence in distribution (see, for example, Corollary 3.3 in chapter 3 of \cite{ek}), so if $\nu$ denotes the distribution of $(R(t_1), \dots, R(t_k))$, then $\rho(\nu_n, \nu) \rightarrow 0$ as $n \rightarrow \infty$ by (\ref{convprob2}).  Combining these results gives
\begin{equation}\label{convprob3}
(d_n |S_{t_1^*}|, \dots, d_n |S_{t_k^*}|) \rightarrow_d (R(t_1), \dots, R(t_k))
\end{equation}
By writing $Z_n(t_i) = (Z_n(t_i) - d_n |S_{t_i^*}|) + d_n |S_{t_i^*}|$ for $i = 1, \dots k$, we can combine (\ref{convprob1}) and (\ref{convprob3}) with Slutsky's Theorem (see Theorem 8.6.1 in \cite{resnick}) to obtain (\ref{fddeq}).
\end{proof}

\section{Relative compactness}
\label{tightsec}

To complete the proof of Theorem \ref{Rayleighthm}, it remains to show that the sequence of processes $(Z_n)_{n=1}^{\infty}$ is relatively compact.  The sequence of processes $(Z_n)_{n=1}^{\infty}$ is relatively compact if the following two conditions hold (see Corollary 7.4 in chapter 3 of \cite{ek}):
\begin{itemize}
\item For all $\epsilon > 0$ and $t \geq 0$, there is a compact set $K$ such that $$\liminf_{n \rightarrow \infty} P(Z_n(t) \in K) \geq 1 - \epsilon.$$

\item For all $\epsilon > 0$ and $T > 0$, there is a $\theta > 0$ such that $$\limsup_{n \rightarrow \infty} P(w(Z_n, \theta, T) \geq \epsilon) \leq \epsilon,$$ where
\begin{equation}
w(Z_n, \theta, T) = \inf_{\{t_i\}} \max_i  \sup_{t,u \in [t_{i-1}, t_i)} |Z_n(t) - Z_n(u)|
\label{wdef}
\end{equation}
and $\{t_i\}$ ranges over all sequences $(t_i)_{i=0}^m$ with $m \geq 1$ such that $0 = t_0 < \dots < t_{m-1} < T \leq t_m$ and $\min_i (t_i - t_{i-1}) \geq \theta$.
\end{itemize}
The first condition follows immediately from the convergence in distribution of $Z_n(t)$ to $R(t)$.  Therefore, our goal in this section is to establish the second condition.

For $t \geq 0$, define
\begin{equation}
g(t) = \left\{
\begin{array}{ll} 
\lfloor a_n |G_n|^{1/2} t \rfloor & \mbox{ in case 1 } \\
\lfloor a_n n^2 (\log n)^{1/2} t \rfloor & \mbox{ in case 2}.
\end{array} \right.
\nonumber
\end{equation}
Also, let
\begin{equation}
L_n = \left\{
\begin{array}{ll} 
a_n |G_n|^{1/2} & \mbox{ in case 1 } \\
a_n n^{2} (\log n)^{1/2} & \mbox{ in case 2}
\end{array} \right.
\nonumber
\end{equation}
and 
\begin{equation}
M_n = \left\{
\begin{array}{ll} 
b_n^{-1} |G_n|^{1/2} & \mbox{ in case 1 } \\
b_n^{-1} n^{2} (\log n)^{1/6} & \mbox{ in case 2}.
\end{array} \right.
\nonumber
\end{equation}
We now choose the points $(t_i)_{i=0}^m$.  Suppose $T$, $\theta$, and $\epsilon$ are fixed.  Let $t_i' = 0$.  For positive integers $i$, let
\begin{equation}
t'_i = \inf \big\{t \geq t'_{i-1} + \theta: X_{g(t)} = X_{g(t) - v} \in LE \big( (X_u)_{u=0}^{g(t) - 1} \big) \mbox{ for some }v \geq M_n \epsilon \big\}.
\label{deftiprime}
\end{equation}
This means that at the time $g(t'_i)$, the random walk $(X_t)_{t \geq 0}$ completes a long loop, causing a downward jump in the process $(Z_n(t), t \geq 0)$ at time $t'_i$.  Since $t'_{i} \geq t'_{i-1} + \theta$ for all $i$, we can choose the points $(t_i)_{i=0}^m$ to coincide with the $t'_i$, and add additional points in the gaps between the $t'_i$ when the gaps have length greater than $2 \theta$.  More precisely, it is possible to choose $0 = t_0 < \dots < t_{m-1} < T \leq t_m$ such that $\theta \leq t_i - t_{i-1} \leq 2 \theta$ for $i = 1, \dots, m$ and, if $t'_j \leq T$, then $t'_j = t_i$ for some $i$.  Note that with this construction, if $t_i$ does not equal $t'_j$ for some $j$, then the open interval $(t_i - \theta, t_{i+1})$ can not contain any of the $t'_j$.

We will use the sequence $(t_i)_{i=0}^m$ to upper bound the right-hand side of (\ref{wdef}).  If $t < u$, then we need both an upper bound for $Z_n(u) - Z_n(t)$, which will show that the loop-erased path does not grow too quickly, and an upper bound for $Z_n(t) - Z_n(u)$, which will show that the loop-erased path does not decrease in length too quickly.  The bounds for $Z_n(u) - Z_n(t)$, provided in Lemmas \ref{upper1} and \ref{upper2}, are relatively straightforward.  Obtaining an upper bound for $Z_n(t) - Z_n(u)$ is more difficult because the length of the loop-erased random walk jumps down at the times when the random walk makes long loops.  However, because the points $t_i'$ are chosen to coincide with these jumps as long as they are far enough apart, the right-hand side of (\ref{wdef}) will only be large if the process $Z_n$ makes two jumps within a time interval of length $\theta$.  This will typically happen only if the random walk makes two nested long loops within a short time, as we show in Lemma \ref{nested}, and we will bound the probability of getting two long loops within a short time in Lemma \ref{Alem}.  Note that Lemmas \ref{upper1}, \ref{nested}, and \ref{Alem} imply Property 2 in case 1, while Lemmas \ref{upper2}, \ref{nested}, and \ref{Alem} imply Property 2 in case 2.  Theorem \ref{Rayleighthm} follows.

\subsection{Bounding the growth of the loop-erased walk}

For case 1, the necessary bound on the growth of the loop-erased random walk follows easily from the fact that the length of the loop-erased walk can grow by at most one on each time step.

\begin{Lemma}
In case 1, there is a constant $C$ such that if $\epsilon > 0$ and $\theta < \epsilon/C$, and if $t_{i-1} \leq t < u < t_i$, then $Z_n(u) - Z_n(t) \leq \epsilon$.
\label{upper1}
\end{Lemma}

\begin{proof}
If $t_{i-1} \leq t < u < t_i$, then $u - t \leq 2 \theta$ and $Y_{g(u)} - Y_{g(t)} \leq g(u) - g(t) \leq 1 + 2 a_n |G_n|^{1/2} \theta$.  Therefore, $$Z_n(u) - Z_n(t) = b_n |G_n|^{-1/2} (Y_{g(u)} - Y_{g(t)}) \leq b_n |G_n|^{-1/2} (1 + 2 a_n |G_n|^{1/2} \theta) \leq C \theta$$ for some constant $C$, which implies the lemma.
\end{proof}

For case 2, this trivial bound is insufficient, and we must make use of the fact that the length of the loop-erasure of a random walk segment of length $r$ is typically of order $r/(\log r)^{1/3}$.  Recall that $A_i = \{(i-1)r, \dots, ir-1\}$.   Let $\ell = \lceil g(T + 2\theta)/r \rceil$.  Let $\eta > 0$, and let $w = \lfloor n^2 (\log n)^{\eta} \rfloor$.  It follows from the proof of Proposition 4.1 of \cite{schw} that for sufficiently large $n$, outside of an event of probability at most $C (\log n)^{-1/11}$, the random walk $(X_t)_{t=0}^{g(T + 2 \theta)}$ satisfies the following two properties.  Recall that $\tau_n$ is the uniform mixing time of the random walk, defined in (\ref{mixdef}).
\begin{itemize}
\item Every interval $[t, t + w]$ with $0 \leq t \leq g(T + 2 \theta) - w$ contains a local cutpoint, that is, a point $u$ such that $\{X_{u - 2 \tau_n}, \dots, X_{u - 1}\} \cap \{X_{u+1}, \dots, X_{u + 2 \tau_n}\} = \emptyset$.

\item For all $i \leq \ell$, if $s, t \in A_i$ with $s + 2 \tau_n \leq t$, then $X_s \neq X_t$.
\end{itemize}

To define the time indices retained after loop-erasure, for positive integers $u \leq v$, let $\sigma_0^{u,v} = \max\{t \in [u, v]: X_t = X_u\}$.  For $i \geq 1$, if $\sigma^{u,v}_{i-1} < v$, let $\sigma^{u,v}_i = \max\{t \in [\sigma_{i-1}^{u,v}, v]: X_t = X_{\sigma_{i-1}^{u,v} + 1}\}$.  Let $W(u,v)$ be the set of times $\sigma_i^{u,v}$, so the path $LE((X_t)_{t=u}^v)$ consists of the points $(X_t)_{t \in W(u,v)}$.  We then have the following result.

\begin{Lemma}
Suppose the two conditions above hold.  Then for all $t \leq g(T + 2 \theta)$ and all $j \leq \ell$, we have $|W(0, t) \cap A_j| \leq |LE((X_t)_{t \in A_j})| + w + 2 \tau_n$.
\label{erasebound}
\end{Lemma}

\begin{proof}
We use ideas from the proofs of Lemma 3.9 and Proposition 3.10 in \cite{schw}.
Note that $|W(0, t) \cap A_j| = 0$ for $t < (j-1)r$, and $t \mapsto |W(0, t) \cap A_j|$ is nonincreasing for $t \geq jr - 1$ because after time $jr - 1$, indices in $A_j$ can no longer be added but can be lost due to loop erasure.  Therefore, we may assume that $(j-1)r \leq t \leq jr - 1$.  Let $$z = \min\{u \geq (j-1)r: X_v \neq X_s \mbox{ for all }v \in [u, t], s \in \{0, 1, \dots, (j-1)r - 1\} \cap W(0, u-1)\},$$ which is the first time during the $j$th segment of length $r$ after which there are no more intersections involving earlier segments.  Then $W(0, t) \cap A_j = W(z-1, t)$ because the indices $(j-1)r, \dots, z-2$ get erased at time $z-1$.  By the first property above, if $z-1 \leq g(T + 2 \theta) - w$, the interval $[z-1, z-1 + w]$ has a local cutpoint, which we call $u$.  We also have $$W((j-1)r, t) \cap \{u, \dots, t\} = W(z-1, t) \cap \{u, \dots, t\}$$ because between times $u$ and $t$, there are no loops involving points before time $u$, as the first property prohibits short loops and the second property prohibits loops of length longer than $2 \tau_n$.
Therefore,
\begin{equation}\label{tech1}
|W(0, t) \cap A_j| \leq |W((j-1)r, t) \cap \{u, \dots, t\}| + w,
\end{equation}
as any discrepancy between $W(z-1, t)$ and $W((j-1)r, t)$ must come from the interval $[z-1, u-1]$, which has length at most $w$.  It follows that
\begin{align} \label{tech2}
|W((j-1)r, t) \cap \{u, \dots, t\}| &\leq |W((j-1)r, t)| \nonumber \\
&\leq |W((j-1)r, jr - 1)| + 2 \tau_n \nonumber \\
&= |LE((X_t)_{t \in A_j})| + 2 \tau_n,
\end{align}
where the second inequality holds because, by the second property above, no point before time $t - 2 \tau_n$ can be erased between times $t$ and $jr - 1$, while the equality is just the definition of $w$.  The result follows from (\ref{tech1}) and (\ref{tech2}).
\end{proof}

\begin{Lemma}
In case 2, there is a constant $C$ such that if $\epsilon > 0$ and $\theta < \epsilon/C$, then for sufficiently large $n$,
$$P \big( \max_i \sup_{t_{i-1} \leq t < u < t_i} (Z_n(u) - Z_n(t)) \geq \epsilon \big) \leq \epsilon.$$
\label{upper2}
\end{Lemma}

\begin{proof}
By (\ref{lenbound}) and Lemma \ref{erasebound}, and the fact that $\tau_n$ is $O(n^2)$ by (\ref{taun2}), there are constants $C_1$ and $C_2$ such that with probability at least $1 - C_1 (\log n)^{1/11}$, we have $|W(0,t) \cap A_j| \leq C_2 n^2 (\log n)^{5/66}$ for all $t \leq g(T + 2 \theta)$ and all $j \leq \lceil g(T + 2 \theta)/r \rceil$.
Suppose $t_{i-1} \leq t < u < t_i$ for some $i$.  We have $$Z_n(u) - Z_n(t) = M_n^{-1} (Y_{g(u)} - Y_{g(t)}) \leq M_n^{-1} (|W(0, g(u)) \cap \{g(t) + 1, \dots, g(u)\}|).$$  
Now $u - t \leq 2 \theta$, so $g(u) - g(t) \leq 1 + 2 a_n n^2 (\log n)^{1/2} \theta$.  It follows that there is a constant $C_3$ such that the number of segments $A_j$ of length $r = \lfloor n^2 (\log n)^{9/22} \rfloor$ that intersect $[g(t) + 1, g(u)]$ is at most $C_3 (\log n)^{1/11} \theta$. 
Therefore, with probability at least $1 - C_1 (\log n)^{-1/11}$, we have $$Z_n(u) - Z_n(t) \leq M_n^{-1} \cdot C_3 (\log n)^{1/11} \theta \cdot C_2 n^2 (\log n)^{5/66} = C_2 C_3 b_n \theta.$$  The result follows by choosing $C > C_2 C_3 \sup_n b_n$ and $n$ large enough that $C_1 (\log n)^{-1/11} < \epsilon$.
\end{proof}

\subsection{Bounding the probability of nearby long loops}

We now work towards bounding the probability that $Z_n(t) - Z_n(u) \geq \epsilon$ if $t_{i-1} \leq t < u < t_i$ for some $i \leq m$.  The next two lemmas show that we typically only have $Z_n(t) - Z_n(u) \geq \epsilon$ if the random walk makes two long loops within a short time.

\begin{Lemma}
If $t_{i-1} \leq t < u < t_i$ and $Z_n(t) - Z_n(u) \geq \epsilon$, then $t_{i-1} = t'_j$ for some $j$.
\label{titjlem}
\end{Lemma}

\begin{proof}
If $Z_n(t) - Z_n(u) \geq \epsilon$, then $Y_{g(t)} - Y_{g(u)} \geq M_n \epsilon$.  This is only possible if a portion of the path at time $g(t)$ of length at least $M_n \epsilon$ is erased by time $g(u)$.  This means that there are integers $v_1$ and $v_2$ satisfying $v_1 \leq g(t) - M_n \epsilon$ and
$g(t) < v_2 \leq g(u)$ such that $X_{v_2} = X_{v_1} \in LE((X_s)_{s=0}^{v_2-1})$.  Let $x = \min\{y: g(y) = v_2\}$.  By (\ref{deftiprime}), either $x = t_j'$ or $t_j' < x < t_j' + \theta$ for some $j$.

Proceeding by contradiction, suppose we do not have $t_{i-1} = t'_j$ for some $j$.  Then, as previously observed, the interval $(t_{i-1} - \theta, t_i)$ does not contain any of the $t_j'$.  However, this contradicts the facts that $t_j' \leq x < t_j' + \theta$ and $t_{i-1} < x < t_i$.
\end{proof}

\begin{Lemma}
Let $\epsilon > 0$.  Let $A$ be the event that there exist integers $v_1, v_2, v_3, v_4$ such that $X_{v_1} = X_{v_4}$, $X_{v_2} = X_{v_3}$, and the following hold:

\smallskip
(a) $v_2 \geq v_1 + \frac{1}{2} M_n \epsilon$.

(b) $v_3 \geq v_2 + M_n \epsilon$.

(c) $v_3 < v_4 \leq v_3 + 1 + 2 L_n \theta$.

(d) $v_4 \leq g(T + 2 \theta)$.
\label{nested}

\smallskip
\noindent In cases 1 and 2, there is a constant $C$ such that if $\theta < \epsilon/C$, then for sufficiently large $n$, $$P \big( \big\{ \max_i \sup_{t_{i-1} \leq t < u < t_i} (Z_n(t) - Z_n(u)) \geq \epsilon \big\} \cap A^c \big) \leq \epsilon.$$
\end{Lemma}

\begin{proof}
Suppose $t_{i-1} \leq t < u < t_i$ for some $i \leq k$, and suppose $Z_n(t) - Z_n(u) \geq \epsilon$.  Suppose also that
\begin{equation}
\max_i \sup_{t_{i-1} \leq x < y < t_i} (Z_n(y) - Z_n(x)) \leq \frac{\epsilon}{2}.
\label{maxbound}
\end{equation}
By Lemmas \ref{upper1} and \ref{upper2}, it suffices to show that under these assumptions, the event $A$ occurs.

By Lemma \ref{titjlem}, there is a $j$ such that $t_{i-1} = t'_j$.  Set $v_3 = g(t'_j)$.  By (\ref{deftiprime}), there is an integer $v_2$ such that (b) holds and $X_{v_3} = X_{v_2} \in LE((X_s)_{s=0}^{v_3-1})$.  This means that the portion of the walk between times $v_2$ and $v_3$ is erased when the loop forms at time $v_3$.  Since $Z_n(t) - Z_n(u) \geq \epsilon$, we have $Y_{g(t)} - Y_{g(u)} \geq M_n \epsilon$.  Also, (\ref{maxbound}) gives $Y_{g(t)} - Y_{v_3} = M_n(Z_n(t) - Z_n(t_j')) \leq \frac{1}{2} M_n \epsilon$.
Therefore, $$Y_{g(u)} = Y_{v_3} + (Y_{g(t)} - Y_{v_3}) + (Y_{g(u)} - Y_{g(t)}) \leq Y_{v_3} + \frac{M_n \epsilon}{2}  - M_n \epsilon \leq Y_{v_3} - \frac{M_n \epsilon}{2}.$$
It follows that between times $v_3$ and $g(u)$, a portion of the walk of length at least $\frac{1}{2} M_n \epsilon$ must get erased.  Since the portion of the walk between times $v_2$ and $v_3$ was already erased at time $v_3$, this can only happen if there is some time $v_4$, with $v_3 < v_4 \leq g(u)$, such that $X_{v_4} = X_{v_1} \in LE((X_w)_{w=0}^{v_4-1})$ and $v_1 \leq v_2 - \frac{1}{2} M_n \epsilon$.  Thus, (a) holds.  Also, (c) holds because $v_4 - v_3 \leq g(u) - v_3 \leq g(t_i) - g(t_{i-1}) \leq 1 + 2 L_n \theta$.  Finally, (d) holds because $v_4 \leq g(t) \leq g(t_k) \leq g(T + 2 \theta)$.  We conclude that $A$ occurs.
\end{proof}

\begin{Lemma}
Let $\epsilon > 0$, and define the event $A$ as in Lemma \ref{nested}.  In cases 1 and 2, there is a constant $C$ such that if $\theta < \epsilon/C$, then $P(A) \leq \epsilon$ for sufficiently large $n$.
\label{Alem}
\end{Lemma}

\begin{proof}
Let $\delta_n = \lfloor \frac{1}{5} M_n \epsilon \rfloor$.  For positive integers $k$ such that $(k-1) \delta_n \leq g(T + 2 \theta)$, let $$B_k = \{(k-1) \delta_n, (k-1) \delta_n + 1, \dots, k \delta_n - 1\}.$$ Also, for such $k$, and for all nonnegative integers $\ell$ such that $(k+1) \delta_n + \ell L_n \theta \leq g(T + 2 \theta)$, let $$D_{k, \ell} = \{(k+1) \delta_n + \lfloor \ell L_n \theta \rfloor, (k+1) \delta_n + \lfloor \ell L_n \theta \rfloor + 1, \dots, (k+1) \delta_n + \lfloor (\ell + 3) L_n \theta \rfloor + 1\}.$$  Note that if $A$ occurs, there must exist integers $j$, $k$, and $\ell$ with $k \geq j + 2$ and $\ell \geq 0$ such that $v_1 \in B_j$, $v_2 \in B_k$, $v_3 \in D_{k, \ell}$, $v_4 \in D_{k, \ell}$, $X_{v_1} = X_{v_4}$, and $X_{v_2} = X_{v_3}$.  In particular, $(X_t)_{t \in D_{k, \ell}}$ intersects both $(X_t)_{t \in B_j}$ and $(X_t)_{t \in B_k}$.

The number of positive integers $k$ such that $(k-1) \delta_n \leq g(T + 2 \theta)$ is at most $1 + g(T + 2 \theta)/\delta_n$, and the number of nonnegative integers $\ell$ such that $\ell L_n \theta \leq g(T + 2 \theta)$ is at most $1 + g(T + 2 \theta)/L_n \theta$.  Therefore, there is a constant $C_1$ such that the number of triples $(j, k, \ell)$ that we must consider is at most $$\bigg( 1 + \frac{g(T + 2 \theta)}{\delta_n} \bigg)^2 \bigg( 1 + \frac{g(T + 2 \theta)}{L_n \theta} \bigg) \leq C_1 \bigg( \frac{L_n}{M_n \epsilon} \bigg)^2 \bigg( \frac{L_n}{L_n \theta} \bigg) = \frac{C_1 L_n^2}{M_n^2 \epsilon^2 \theta}.$$

In cases 1 and 2, when $k \geq j + 2$ and $\ell \geq 0$, for sufficiently large $n$ there is a gap between $B_k$ and $B_j$ of length greater than $\tau_n$ and a gap between $B_k$ and $D_{k, \ell}$ of length greater than $\tau_n$.  Therefore, in case 1, it follows from Lemma 5.2 and equation (8) of \cite{perrev} that there is a constant $C_2$ such that the probability that $(X_t)_{t \in D_{k, \ell}}$ intersects both $(X_t)_{t \in B_j}$ and $(X_t)_{t \in B_k}$ is at most $C_2 (L_n \theta)^2 (M_n \epsilon)^2/|G_n|^2$.  Likewise, in case 2, it follows from Proposition 3.3 of \cite{schw} that there is a constant $C_3$ such that 
the probability that $(X_t)_{t \in D_{k, \ell}}$ intersects both $(X_t)_{t \in B_j}$ and $(X_t)_{t \in B_k}$ is at most $C_3 (L_n \theta)^2 (M_n \epsilon)^2/(n^8 (\log n)^2)$.  Putting together these bounds, we see that in case 1, we have $$P(A) \leq C_1 C_2 \bigg( \frac{L_n^2}{M_n^2 \epsilon^2 \theta} \bigg) \bigg( \frac{L_n^2 M_n^2 \theta^2 \epsilon^2}{|G_n|^2} \bigg) = \frac{C_1 C_2 \theta L_n^4}{|G_n|^2} \leq C_4 \theta$$ for some constant $C_4$.  Likewise, in case 2, we have
$$P(A) \leq C_1 C_3 \bigg( \frac{L_n^2}{M_n^2 \epsilon^2 \theta} \bigg) \bigg( \frac{L_n^2 M_n^2 \theta^2 \epsilon^2}{n^8 (\log n)^2} \bigg) = \frac{C_1 C_3 \theta L_n^4}{n^8 (\log n)^2} \leq C_5 \theta$$ for some constant $C_5$.  These bounds imply the lemma.
\end{proof}

\section*{Acknowledgments}
\noindent The author thanks Yuval Peres, Jim Pitman, and David Revelle for helpful discussions.  He also thanks the referee for suggestions which improved the exposition of the paper.

\end{document}